\documentclass[twoside, 10pt]{article}
\usepackage{latexsym}
\def \Chi{\ensuremath{\mathbf{\chi}}}
 
\def \H{\ensuremath{\Bbb H}} 
\def \BB{\ensuremath{\mathcal B}} 
\def \RR{\ensuremath{\mathcal R}} 
\def \VV{\ensuremath{\mathcal V}} 
 
\def \NN{\ensuremath{\mathcal N}} 
\def \LL{\ensuremath{\mathcal L}} 
\def \DD{\ensuremath{\mathcal D}} 
\def \ZZ{\ensuremath{\mathcal Z}} 
\def \CC{\ensuremath{\mathcal C}} 
\def \UU{\ensuremath{\mathcal U}}
\newcommand {\HP}[1][n]{\ensuremath{\Bbb {HP}^{#1}}} 
 
\newcommand {\Tn}[2][\id]{\ensuremath{\mathrm{T}_{#1}{#2}}} 
\newcommand {\CTn}[2][\id]{\ensuremath{\mathrm{T}^{\ast}_{#1}{#2}}}

\newcommand \Sp[1][n]{\ensuremath{\mathrm{\mathbf{Sp}}(#1)}} 
\newcommand \SL[1][n]{\ensuremath{\mathrm{\mathbf{SL}}(#1)}} 
\renewcommand \sp[1][n]{\ensuremath{\frak{sp}_#1}} 
\def \hh {\ensuremath{\frak{h}}}
\def \Sham{{\aleph}}
\def \Sph{\ensuremath{\Sigma}} 
\def \sph{\ensuremath{\frak{s}}}
\def \s{\ensuremath{\frak{s}}}
 
\def \R{\ensuremath{\Bbb R}} 
\def \Z{\ensuremath{\Bbb Z}} 
\def \sham{\ensuremath{\kappa}} 
\def \Wedge{\ensuremath{\bigwedge}}

\renewcommand {\dim}[2][\R]{\ensuremath{\mathrm{\mathbf{dim}}_{#1}(#2)}}
\newcommand {\sgn}[1]{\ensuremath{\mathrm{\mathbf{sgn}}(#1)}}

\newcommand {\diag}[1]{\ensuremath{\mathrm{\mathbf{diag}}(#1)}}
\newcommand {\len}[1]{\ensuremath{\mathrm{\mathbf{len}}(#1)}}
\renewcommand {\det}[1]{\ensuremath{\mathrm{\mathbf{det}}(#1)}}
\renewcommand {\exp}[1]{\ensuremath{\mathrm{\mathbf{exp}}(#1)}}
\newcommand {\ad}[2]{\ensuremath{\mathrm{\mathbf{ad}}_{#1}(#2)}}
\def \id{\ensuremath{\mathbf{e}}}
\newcommand {\GL}[2][n]{\ensuremath{\mathrm{\mathbf{GL}}_{#1}(#2)}} 
\newcommand {\gl}[2][n]{\ensuremath{{\frak{gl}}_{#1}(#2)}} 
\newcommand {\mat}[4]{\left(\begin{tabular}{c c} \ensuremath{#1} & \ensuremath{#2} \\
						 \ensuremath{#3} & \ensuremath{#4} \\
			    \end{tabular}\right)}
\def \i{{\ensuremath{\bf{i}}}}
\def \j{{\ensuremath{\bf{j}}}}
\def \k{{\ensuremath{\bf{k}}}}
 
\def \proof{{\noindent{\it Proof.\ \ }}} 
\def \example{{\noindent{\it Example.\ \ }}} 
\def \qed{$\Box$}
\def \qex{$\bigcirc$}
\newtheorem{Th}{Theorem}[section] 
\newtheorem{definition}[Th]{Definition} 
\newtheorem{prop}[Th]{Proposition} 
\newtheorem{lemma}[Th]{Lemma} 
\newtheorem{lem}[Th]{Lemma}

\date{March 28, 2002}
 
\title{Geometry of Four-vector Fields on \\
Quaternionic Flag Manifolds.} 

\author{Philip Foth \ \ \& \  
Frederick Leitner}

\begin{document} 
 
\maketitle 
\input amssym.def 
\markboth{Foth, Leitner}{Quatrisson structures and $\H$-flags} 
\setcounter{equation}{0}

\begin{abstract}  The purpose of this paper is to describe certain natural
$4$-vector fields on quaternionic flag manifolds, which 
geometrically determine the Bruhat cell decomposition. These structures
naturally descend from the symplectic group \Sp, and are 
related to the dressing action given by the Iwasawa
decomposition of the general linear group over the quaternions, $\GL{\H}$. 

\end{abstract}

%\footnotetext{\noindent{\it AMS subj. class.}: \ \ 
%primary 58F05, secondary 57R25, 22E30.} 

\section{Introduction}

In this paper we wish to describe certain natural $4$-vector fields 
on quaternionic flag manifolds. In the context of the 
Poisson geometry,  a bi-vector field is penultimate in the
study of the geometry of the underlying manifold.  
Analogously, we make use of a
$4$-vector field,
closed under the Schouten bracket with itself, 
 which we call a  \emph{quatrisson} structure, to 
reveal the internal structures
of certain natural spaces arising in geometry, namely
quaternionic flag manifolds. 
A more general definition involving a multi-vector field was first 
given in \cite{Perelomov}. Another generalization, the Nambu-Poisson structure,
was studied in \cite{Takh}. 

Quaternionic flag manifolds possess
natural group invariant quatrisson 
structures,
and the study of the geometry of the flag manifolds can be pursued in the natural
setup of quatrisson 4-vector fields and tetraplectic structures \cite{F1}. 
In particular, we describe 
the so-called Bruhat quatrisson 4-vector fields on quaternionic 
flag manifolds where the leaf decompositions
coincide with the Bruhat decompositions 
of  $\GL{\H}$ defined purely combinatorially. 
We also show that the existence of the Bruhat decomposition
leads to a description of the tetraplectic leaves in the group $\Sp$ 
in terms of the dressing action on the group. 

Drinfeld \cite{Dr}, Lu and Weinstein \cite{LW}, 
Semenov-Tian-Shansky \cite{STS}, and Soibelman \cite{Soibel}
first described this setup in the context of standard Poisson geometry,
and this viewpoint has been elaborated by many others. 
Several important features of the Poisson geometry of flag manifolds
readily translate to our situation, including
Schubert calculus and a version of generalized hamiltonian dynamics.
We suggest that further
studies of these structures might lead to interesting results related 
to the geometry and (equivariant) differential calculus on quaternionic 
flag manifolds as well as quantum groups. 

\ 

\noindent{\bf Acknowledgements.} The first author is grateful to 
Sam Evens and Lu Jiang-Hua for many conversations related to Poisson
geometry. The first author was supported by NSF grant DMS-0072520.
The second author was supported by an NSF VIGRE graduate fellowship.

\section{Quaternionic matrices and flags} 
We begin with some generalities on quaternionic matrices for which we
define the following subgroups of $\GL{\H}$:
\begin{eqnarray*}
	\RR 
		& 
		:= 
			&  
			\{ \diag{r_1,\dots,r_n} \ | \ r_i \in \R_+ \} 
	\\
	\UU 
		& 
		:= 
			& 
			\{\mbox{upper triangular matrices with  }1\mbox{'s along the
			diagonal}\} 
	\\
	\VV_w 
		&
		:=
			&
			\{ U\in \UU \ | \ \left( P_w U P_w^{-1}\right)^{t} \in \UU \}
	\\
		&
			&
			\mbox{here} \ P_w \
			\mbox{denotes the permutation matrix}
	\\
		&
			&
			(P_w)_{i,j}=\delta_{i,w(j)} \ \mbox{for} \
			w\in S_{n}  
	\\ 
        \DD 
		& 
		:= 
			& 
			\{ \diag{d_1,\dots, d_n} \ | \ d_i \in \H^* \}
	\\
        \BB 
		& 
		:= 
			& 
			\UU\DD
\end{eqnarray*}

\noindent Now for $G\in\GL{\H}$ we recall
 the \emph{strict Bruhat normal form} \cite{Draxl} of $G$ as:

\begin{displaymath}
	G = U D P_w V
\end{displaymath}

\noindent Here all the matrices are uniquely determined:
the matrix $D=\diag{d_1,\dots,d_n}$ belongs to $\DD$; $P_w$
is, as usual,  the permutation matrix corresponding to $w\in S_n$;
both $U$ and $V$ belong to $\UU$;  
and we further require that $P_w V P_w^{-1}$ 
is lower triangular with $1$'s along the diagonal, i.e. 
$V\in \VV_w$.  
This decomposition allows us to  define the Dieudonn\'{e} 
determinant \cite{Dieu} as the residue
of $\sgn{w}\cdot\prod d_i$ in $\H/[\H,\H]=\R_{+}$.
Moreover, by means of the strict Bruhat normal form, we obtain the Bruhat decomposition: 
 
\begin{eqnarray*} 
\GL{\H}=\coprod_{w\in S_n} \BB P_w \VV_w. 
\end{eqnarray*}

\noindent Denoting $\ZZ_{w}:= \{ \BB P_w \VV_w \}$, we obtain from the
Bruhat decomposition a parameterization 
of $\GL{\H}$ by $S_n$. 
The condition that  $V^{w }:=wVw^{-1}$ is lower triangular 
implies in the case of $w =\id$ that  $V^{\id}=V$ must be both upper and lower
 triangular  hence equals  the identity matrix and hence
$\dim{\ZZ_{\id}}=\dim{\BB}$.
Taking $w$ to be  the longest permutation, $w_l=\left(n\ (n-1) \cdots \ 2\ 1\right)$,
rotates the matrix $V$ by $180^\circ$ so that it is lower triangular.
As no further conditions
on the entries of $V$ are imposed, we have that
$\dim{\ZZ_{w_l}}=4n^2$.  In general, we
define the length  of a permutation $w$, $\len{w}$,
to be the minimal number of adjacent transpositions required
in  a factorization of the permutation.  
One readily sees that the maximal number of non-zero
entries allowed in $V$, so that $V^w $ is lower triangular,
 is exactly $\len{w }$ so that
$\dim{\ZZ_w }=4\cdot \len{w }+\dim{\BB}$.  
We will see later that
the entries of $V$ give coordinates on the Bruhat cells of quaternionic flags.
 
Denoting conjugate transpose by $(\cdot)^\ast$, we have the Lie group:
\begin{displaymath}
  	\Sp:=\{g\in\GL{\H} \ | \ g^{\ast}g=\id\}
\end{displaymath}

\noindent We also identify the corresponding Lie algebra:

\begin{displaymath}
  	\sp:=\{X\in\gl{\H} \ | \  X + X^{\ast}=0\}
\end{displaymath}
with  $\Tn[x]{\Sp}$ for $x\in\GL{\H}$ by
left translation.  
One knows that matrices in $\Sp$  have Dieudonn\'{e} determinant 1, and
thus lie inside the semi-simple group:
\begin{displaymath}
	\SL:=\{g\in\GL{\H} \ | \ \det{g} = 1\}.
\end{displaymath}

\noindent We define the \emph{spheroid} to be 
$\Sph:=\Sp[1]^{n}\simeq (S^3)^n$, whose elements
are of the form $\diag{\exp{s_1},\dots,\exp{s_n}}$ where the $s_i$ 
are purely quaternionic (no real component).  
This is a subgroup of $\DD$, and we have, in fact, that:
\begin{displaymath}
	\DD = \Sph\times \RR
\end{displaymath}

\noindent We  denote the corresponding Lie algebra by $\s$.
The full quaternionic flag of $\H^n$, which we denote $F_n$,
can now be identified as $F_n \simeq \Sph\backslash\Sp\simeq
\BB\backslash\GL{\H}$.
Using the second identification of $F_n$, we denote
$\CC_w := \BB\backslash \ZZ_w$, which we call the \emph{Bruhat cells} of the
flag.  By our discussion above, we see that
\ $\dim{\CC_w}= 4\cdot \len{w}$. 

\
 
\example
Consider the space $\HP[1]$ identified as
 $\HP[1]\simeq \Sph\backslash\Sp[2] \simeq S^{4}$
from which we  obtain the  fibration: 
$$
\begin{array}{ccc}
\Sph\simeq S^3\times S^3\simeq\Sp[1]\times\Sp[1]	& \rightarrow	& \Sp[2]  \\
	             		&             	&  \downarrow      \\
	       			&             	& \HP[1] \simeq S^4
\end{array}.
$$

The Bruhat decomposition yeilds a decomposition of 
$\Sph\backslash\Sp[2]\simeq\HP[1]$ into the cells
$\CC _{(12)}$ and $\CC_\id$
which have real dimensions $4$ and $0$ respectively.
We view the cells under the identifications that
	$\CC_{(12)}\simeq \NN^2 = \mat{1}{\ast}{0}{1}\simeq\H$ 
and $\CC_\id$ is the North pole.  

Recall that $S^{4}$ has neither
a symplectic nor a complex structure (nor even an almost complex \cite{Mass}).  
This is one of the reasons for introducing the 
tetraplectic structure in \cite{F1}. 
\qex

\section{ Quatrisson and tetraplectic structures}

Recall that a symplectic manifold is a manifold equipped with a closed 
non-degenerate 2-form. We also recall 
that a Poisson manifold is a 
manifold equipped with a  bi-vector field that induces a Lie algebra structure 
on the space of smooth functions,  compatible with the commutative
product of functions via the Leibniz rule.
In the case of quaternionic flags, 
we make use of the following structures which reflect the  underlying geometry:

\begin{definition} \cite{F1}
	Let $X$ be a real orientable manifold of dimension $4m$.  
	A {\emph tetraplectic} structure on X is a  four-form, $\psi$ satisfying: 

	\begin{tabular}{c l } 
		$1)$ & $\psi$ is closed ($d\psi=0$)\\ 
		$2)$ & $\psi^{m}$ is a volume form \\ 
	\end{tabular} 
	
	\noindent We call the pair
	 $(X, \psi)$  a \emph{tetraplectic manifold}. A map
	$\phi: (X, \psi)\to (X', \psi')$ is called \emph{tetraplectic} if 
	$\phi^*\psi'=\psi$. If, in addition, $\phi$ is a diffeomorphism, then
	we call $\phi$ a \emph{tetraplectomorphism}. 
\end{definition}

\example 
	Let $\psi$ be an $\Sp[2]$-invariant volume form on $S^{4}$.
	Then ($\HP[1]$,$\psi$)  is 
	a teraplectic manifold.
	In fact, in \cite{F1} the construction of invariant tetraplectic structures
	on all quaternionic flag manifolds was given. 
\qex 

\
 
\noindent One can define a standard Poisson structure 
on a manifold by giving a bi-vector field whose
the Schouten bracket with itself is zero. 
However, in order to reflect the geometry 
of our situation we will make use of 4-vector fields, for
which we recall the following \cite{Vais}:  

\begin{prop} 
	Denoting $\Wedge^i\Chi(M)$ 
	the space of $i$-vector fields on $M$,	
	there exists a unique  bracket, called the \emph{Schouten bracket}:
	\begin{displaymath}
		[\cdot,\cdot]:\Wedge^{p}\Chi(M)\times\Wedge^{q}\Chi(M) 
		\rightarrow\Wedge^{p+q-1}\Chi(M)
	\end{displaymath} 
	which extends the usual Lie bracket of vector fields and is an 
	$\R$-linear operation satisfying the following identities:  

	\begin{tabular}{l c r} 
		$1)$ 
			& 
			$[P,Q] = (-1)^{pq}[Q,P]$ 
				& 
				(Anti-Symmetry) 
		\\
		$2)$ 
			& 
			$[P,Q\wedge R] = [P,Q]\wedge R  + (-1)^{pq+q}Q\wedge[P,R]$ 
				& 
				(Leibniz)
		\\
       	 	$3)$ 
			& 
			$(-1)^{p(r-1)}[P,[Q,R]] + (-1)^{q(p-1)}[Q,[R,P]]$  
				& 
				{} 
		\\ 
	  	{}   
			& 
			$+ (-1)^{r(q-1)}[R,[P,Q]] = 0$  
				& 	
				(Jacobi)
	\end{tabular}  
\end{prop}

\noindent We recall that in \cite{Perelomov} the authors  
use the  vanishing of the Schouten bracket  
of a $p$-vector field $\xi$ with itself, $[\xi,\xi]=0$,
to define the Generalized Poisson Structures (GPS).
For etymologic-semantic reasons, we give 
the following definition:

\begin{definition} 
	Let $M$ be a manifold, and let $\xi$ be a 4-vector field 
	on $M$ satisfying $[\xi, \xi]=0$. We call $\xi$ a \emph{quatrisson} 
	structure on $M$ and the pair $(M, \xi)$ a \emph{quatrisson manifold}. 
\end{definition}

\begin{definition} 
	For two quatrisson manifolds ($X$, $\xi$) and ($X'$,$\xi'$) 
	a map $\phi:X\rightarrow X'$ is called a \emph{quatrisson map} if 
	for any quadruple of functions $f_i\in C^{\infty}(X')$, 
	$1\le i\le 4$
	the following identity holds: 
	\begin{displaymath}
		\xi( d\phi^*f_1\wedge d\phi^* f_2\wedge d\phi^* f_3\wedge d\phi^* f_4)=
		\phi^* \xi'(df_1\wedge df_2\wedge df_3\wedge df_4).
	\end{displaymath}
	i.e. $\phi_*(\xi)=\xi'$. 
\end{definition}

\begin{definition} 
	Let $\xi$ be a  4-vector field  on
	a $4m$-dimensional manifold, $M$.
	Then we call $\xi$  
	\emph{non-degenerate} if 
			$\xi^{\wedge m}$ is a nowhere vanishing 
$4m$-vector field. 
\end{definition} 

If $\xi$ is a non-degenerate vector field on $M$, then 
$\xi$ induces a surjection $\Wedge^3\CTn[x]{M}\rightarrow\Tn[x]{M}$ 
for all $x\in M$,  obtained by contraction with $\xi$.  
If $(M,\xi)$ is quatrisson, we define the \emph{rank} of $\xi$ at $x\in M$ 
as the dimension of the image of this map.
One can see that a quatrisson structure $\xi$ on $M$ is 
non-degenerate if the rank of $\xi$ at any point of $M$ is equal to the
dimension of $M$.

\begin{lem} 
	Letting $(M, \xi)$ be as above, the 
	rank of $\xi$ at any $x\in M$ is divisible by 4. 
\end{lem} 

\proof This is an easy exercise in multi-linear algebra 
for the reader.   
\qed

\begin{definition} 
	Let $M$ be a manifold equipped with a 4-vector field
	$\xi$. We say that a smooth $4l$-dimensional 
	submanifold, $L$, is a \emph{tetraplectic leaf} in $M$ if:
	\begin{tabular}{c l}
		$1)$ 	
			&
		 	$\xi$ comes from  
			$\Wedge^{4}\Chi(L)$ at all points of $L$
		\\
		$2)$
			&
			$\xi$ is non-degenerate on $L$
		\\
		$3)$
			&
			$L$ is not properly included in any other such submanifold of $M$
		\\
		$4)$
			&
			the four-form, $\psi$,  given by $i_{\psi}\xi^{m}=\xi^{m-1}$, defines a 
		\\
		{}	
			&
			tetraplectic structure on $M$. 
	\end{tabular}
\end{definition}

\noindent To each triple of functions ${\bf f}$=($f_1$, $f_2$, $f_3$), 
we can associate a ``hamiltonian'' vector field 
$X_{\bf f}$ given by $\iota(df_1\wedge df_2\wedge df_3)\xi$.
Then we get the \emph{characteristic distribution} of $M$. 
Unlike in the Poisson case, we cannot
expect in general that $(M, \xi)$ is stratified as a union of 
smooth tetraplectic leaves, even if $\xi$ is quatrisson,
see Example 8 in \cite{Ibanez2}. 
However, the particular case of this result 
for quaternionic flag manifolds will follow later.

\

\example Let $F_n$ be a quaternionic flag manifold considered as a 
tetraplectic manifold with an $\Sp[n]$-invariant $4$-form $\psi$ \cite{F1}. 
The corresponding  4-vector field, $\chi$, defined by
 $i_\chi(\psi^m)=\psi^{m-1}$,
 is quatrisson and $\Sp$-invariant. We refer to $\chi$
as the invariant quatrisson structure of $F_n$.
\qex

\section{Quatrisson structures on ${\bf HP}^1$.} 
 
For our flag manifolds, we construct a Bruhat quatrisson structure explicitly by 
analogy to the  Poisson case as in \cite{LW} or \cite{Soibel}. 
We begin by defining a quatrisson 4-vector field, $\sham$, on $\Sp[2]$, which  
we show descends to the quotient $\Sph\backslash\Sp[2]$.   
 
We begin by defining an element $\Lambda$ of $\wedge^4\sp[2]$ 
in terms of the following basis for $\sp[2]$: 
\begin{eqnarray*} 
	E = \mat{0}{1}{-1}{0}	  
			&   	
			S_x = \mat{0}{x}{x}{0}  
	\\ 
	H_x = \mat{x}{0}{0}{-x}	  
			&  	
			M_x = \mat{x}{0}{0}{x}	  
\end{eqnarray*} 
where $x$ is one of $\{ \i, \j, \k\}$.  For convenience we denote
$S_{-x}:=-S_x$, $H_{-x} := -H_x$, and $M_{-x} := - M_x$.  We now can note the
following commutator relations: 
\begin{eqnarray*}
	{[M_x,E]}=0
		&
	{[H_x,E]}=2\cdot S_x
			&
			{[ S_x , E ]}=	2 \cdot H_x
\end{eqnarray*} 
\begin{eqnarray*}
	{[ M_x , M_y ] } = {[ H_x , H_y ]} =  {[ S_x , S_y ] }
		&
		= 
			&
			\Bigg\{ 
			\begin{tabular}{c l}
				$2\cdot M_{x\cdot y} $
					&
					$x\ne y$
				\\
				$0$
					&
					$x=y$
			\end{tabular} 
	\\
	{[S_x,H_y] }
		&
		=
			&
			\Bigg\{ 
			\begin{tabular}{c l}
				$0$
					&
					$x\ne y$
				\\
				$-2\cdot E$
					&
					$x=y$
			\end{tabular}
	\\
	{[S_x,M_y]}
		&
		=
			&
			\Bigg\{ 
			\begin{tabular}{c l}
				$0$
					&
					$x =  y$
				\\
				$2\cdot S_{x\cdot y}$
					&
					$x \ne y$
			\end{tabular}
	\\
		{[H_x,M_y]}
		&
		=
			&
			\Bigg\{ 
			\begin{tabular}{c l}
				$0$
					&
					$x =  y$
				\\
				$2\cdot H_{x\cdot y}$
					&
					$x \ne y$
			\end{tabular}
\end{eqnarray*}

\noindent We may now define: 
\begin{displaymath} 
	\Lambda :=E \wedge S_{\i} \wedge S_{\j}\wedge S_{\k}. 
\end{displaymath} 

\noindent and denoting by $\Lambda^{L}$ and $\Lambda^{R}$  the left 
and right 
invariant $4$-vector fields on $\Sp[2]$ with  value $\Lambda$ 
at the identity element, we let:
\begin{displaymath}	
	\sham=\Lambda^{L}-\Lambda^{R}.
\end{displaymath}
 
\begin{prop} 
	The 4-vector field $\sham$ is a quatrisson structure
	on $\Sp[2]$.  Moreover, $\sham$ 
	descends to  a vector field, $\Sham$,
	on 	$\HP[1]\simeq\Sph\backslash\Sp[2]$  inducing a $\Sph$-invariant quatrisson 
	structure on $\HP[1]$ called the 
	\emph{Bruhat quatrisson structure} \label{prop444} 
\end{prop} 

\proof 
	The fact that $\sham$ is a quatrisson structure on $\Sp[2]$ is
	nothing more than the fact that $[\sham,\sham]\in\Wedge^{7}\Chi(\Sp[2]) = 0$.
	To show that $\Sham$ is $\Sph$-invariant, and hence descends, we may apply
	the same formalism of the Poisson case
	and show that for any $X\in \sph$ we have $\ad{X}{\Lambda} = [X,\Lambda]=0$.
	This follows readily from the above commutator relations and 
	the Leibniz rule of the Schouten bracket.  It is clear that 
	$[\Sham,\Sham]=0$.
\qed

\

\noindent We can make use of the Bruhat decomposition to describe the vector 
field explicitly.  As above, we denote by $\CC_\id$ and $\CC_{(12)}$ 
the cells of $F_n$ corresponding to the 
North pole and the $\H$ components.
It is clear, that at the North pole $\sham$ is the zero vector.  
For $x\in \CC_{(12)}$ we
choose a convenient coset representative in $\Sp[2]$, 
namely:
\begin{displaymath}
	k_x = \frac{1}{\sqrt{1+\rho^2}} \cdot \mat{-\bar{v}}{1}{1}{v},  
\end{displaymath}
where $\rho = |v| = \sqrt{v\bar{v}}$. In fact, the identification of 
$S^4$ as a natural $SO(5)\simeq \Sp[2]/(\Z/2)$ -invariant submanifold of 
$\R^5$ with $\H$ plus the point at infinity using stereographic 
projection sends $v\in \H$ to a point in $S^4$ at the height
$\displaystyle{1-{1\over 1+|v|^2}}$ and the 
same $\Sp[1]$-angular coordinate. 

To compute $\Sham$ at $x$, we identify $\Tn[x]{\Sp[2]}$ with 
$\Tn{\Sp[2]}$ by right translations so that  we have
$\Lambda^{R} = \Lambda$ and  $\Lambda^{L}$ is simply 
conjugation of  $\Lambda$ by $k_x$.  We have thus expressed 
$\Sham=f(v)\partial_{v_1}\wedge\partial_{v_2}\wedge\partial_{v_3}\wedge\partial_{v_4}$
 in terms
of the coordinates $v=v_1+v_2\cdot \i+v_3 \cdot \j + v_4 \cdot \k$.  We would
like more natural coordinates for $\H$, namely if
$v=\rho\cdot\exp{\theta_1\cdot \i}\exp{\theta_2 \cdot \j}\exp{\theta_3 \cdot \k}\in \H$,
then we have $\Sham=$
$\displaystyle{\frac{g(\rho)}{\rho^3}\partial_\rho\wedge\Theta}$, for
 $\Theta=\partial_{\theta_1}\wedge\partial_{\theta_2}\wedge\partial_{\theta_3}$ -
the $\Sp[1]$-invariant $3$-vector field on $\Sp[1]\simeq S^3$. 
To find $g(\rho)$, we divide $\Sham$ by its value at $v=0$, the South pole. 
After a  computer assisted 
computation\footnote{We thank Klaus Lux and Stephane Lafortune for help with this.} 
we see that  
$\displaystyle{g(\rho)=\frac{(1+3\rho^4)}{(1+\rho^2)^3}}$. Thus we have proved that:

\begin{prop} 
	The invariant quatrisson structure on $S^4\simeq \HP[1]$
	is given by: 
	\begin{displaymath}
		\chi := {(1+\rho^2)^4\over \rho^3}\partial_\rho\wedge\Theta
	\end{displaymath}
	and the Bruhat quatrisson structure is given by:
	\begin{displaymath}
		\Sham = {(1+\rho^2)(1+3\rho^4)\over \rho^3}\partial_\rho\wedge\Theta.
	\end{displaymath}
	In particular we have that:
	\begin{displaymath}
	 	\Sham={1+3\rho^4\over (1+\rho^2)^3}\chi.
	\end{displaymath}
\end{prop} 

One can easily see that $\Sham$ has rank four everywhere except for
at the North pole, where it vanishes. Thus the two cells are characterized by 
the rank of $\Sham$.

\section{Quatrisson structures on  flag manifolds} 
 
Following \cite{Lu2} 
we will produce the quatrisson structure on the full flag of 
$\H^{n}$ by way of the so-called multiplication formula.  By analogy to the $\Sp[2]$
 case, for $1\le p<q \le n$ 
 we denote by $E^{p,q}$ the quaternionic matrix whose entries are $0$'s
 everywhere except in the 
($p$,$q$)-th position which is $1$, and the ($q$,$p$)-th position which is $-1$. We also let  
$S_x^{p,q}$ denote the matrix with $0$'s  everywhere except in the 
($p$,$q$)-th and ($q$,$p$)-th positions where the entries are $x$ where $x$ is again
chosen from $\{\i, \j, \k\}$.  Similarly,
these matrices are clearly in the Lie algebra, $\sp$, of $\Sp$, 
and correspond to ``positive roots''
(i.e. pairs of integers $1\le p< q\le n$) as in \cite{Lu2}.   We  define
 $\Lambda\in\Wedge^4\Tn{\Sp}$ by:
\begin{displaymath} 
	\Lambda=\sum_{p<q} E^{p,q}\wedge S_{\i}^{p,q}\wedge S_{\j}^{p,q}\wedge S_{\k}^{p,q}.
\end{displaymath} 
Then if $\Lambda^{L}$ and $\Lambda^{R}$ are the right and left invariant 
$4$-vector fields on $\Sp$ with the values  $\Lambda$ at the identity element
on $\Sp$, we let:
\begin{displaymath}
	 \sham=\Lambda^{L}-\Lambda^{R}.
\end{displaymath}

Unlike the $\Sp[2]$ case, when $n> 2$, one can readily check that $\sham$ will
not be a quatrisson structure on $\Sp$ by making use of the Leibniz rule and commutator
relations similar to those as above and noting that there will be some terms that will not cancel.
However, we still have:

\begin{prop} 
	The 4-vector field $\sham$ descends to $\Sph\backslash\Sp$, inducing a
	\Sph -invariant quatrisson structure, $\Sham$, 
	called the \emph{Bruhat quatrisson structure}.
\end{prop} 

\proof 
	For $\sham$ to descend and be invariant we need to show that both
	the left and right translations by elements of the Spheroid 
	leave $\sham$ invariant, meaning that the 
	adjoint action by the Spheroid on $\Lambda$ is trivial.
	This can be checked similarly to the $n=2$ case  of Proposition
	\ref{prop444}. One can also directly check that $[\Sham, \Sham]=0$ on
	$\Sph\backslash\Sp$, which will also follow from Proposition
\ref{prop:p247}.
\qed

\

\noindent We recall:

\begin{definition}\label{def:mult-act}
	Let $H$ be a Lie group equipped with a 
	multiplicative 4-vector field $\mu$, 
	which acts on a quatrisson manifold $(P, \xi)$:
	\begin{displaymath}
		\beta:H\times P\to P.
	\end{displaymath}
	We say that $H$ acts \emph{multiplicatively} if,
	denoting the corresponding translation maps:
	\begin{eqnarray*}
	 	\beta_h:P\to P
			&
				&
				\beta_y:H\to P
		\\
		y\mapsto h\cdot y
			&
				&
				h \mapsto h\cdot y
	\end{eqnarray*}
	we have:
	\begin{displaymath}
		\xi(h\cdot x) = \beta_{h_{\ast}}\xi(x) + \beta_{x_{\ast}}\mu(h).
	\end{displaymath}
	We sometimes say that the actions is multiplicative with respect to the
	direct sum $4$-vector field $\mu\oplus\xi$ on $H\times P$.
\end{definition}
 
\noindent and notice the following fact (cf. \cite{LW}, \cite{KS}):

\begin{lemma}
	The $4$-vector field $\sham$ on $\Sp$ is 
	multiplicative.
\end{lemma}

\begin{prop} 
	Let $\Sham$ be the Bruhat quatrisson structure on $\Sph\backslash\Sp$.  
	The action map: 
	\begin{displaymath} 
		\Sp \times \Sph\backslash\Sp \rightarrow \Sph\backslash\Sp 
		: (g,h)\mapsto g\cdot h 
	\end{displaymath} 
	is multiplicative with respect to the four-vector field $\sham\oplus\Sham$ 
	on $\Sp\times(\Sph\backslash\Sp)$. 
\end{prop} 
\proof Straightforward.
\qed

\

\noindent We will also make use of the following embeddings:
\begin{eqnarray*}
	f_{r,r+1}:  \Sp[2]  
		&
		\hookrightarrow  
			&
			\Sp 
	\\
	A
		&
		\mapsto 
			&
			A_{r,r+1}
\end{eqnarray*}
where $1\le r <n$ and for  
$A=\mat{a}{b}{c}{d}$  the matrix $A_{r,r+1}$ is given by: 
\begin{displaymath} 
	  \left( \begin{tabular}{c c c c}
		$Id$	 
			&
				&
					&
					$0$
		\\
			&
			$a$
				&
				$b$
					&
		\\
			&
			$c$
				&
				$d$
					&
		\\
		$0$ &
				&
					&
					$Id$
	\end{tabular}\right)
	\begin{tabular}{c}
		$\leftarrow r^{th}$ row 
		\\ 
		\\ 
		\\
	\end{tabular} 
\end{displaymath} 
 
\

\begin{lemma}
	The embeddings $f_{r,r+1}: \Sp[2] \hookrightarrow \Sp$ respect the
	multiplicative 4-vector fields $\sham$. 
\end{lemma}
\proof
	Straightforward. 
\qed

\begin{prop}
	Every tetraplectic leaf $L$ of $\Sp$ lies entirely in some $\ZZ_w$. 
	If $L_w$ is a tetraplectic leaf containing the permutation matrix $P_w$ 
	corresponding to some $w\in S_n$, and we  write
	$w=\prod_{i=1}^{m}\tau_i$ as a minimal product
	of adjacent transpositions, we have a tetraplectomorphism:

	\begin{tabular}{c c l c}
		$F_w:$
			& 
			$L_{\tau_1}\times\dots\times L_{\tau_m}$ 
				& 
				$\rightarrow$ 
					& $L_w$ 
		\\
	           	& 
			$(l_1,\dots,l_m)$ 
				& 
				$\mapsto$     
					& 
					$l_1 l_2 \dots l_m$ 
	\end{tabular}

	\noindent Moreover,  for $\sigma\in\Sph$, 
	the tetraplectic leaf through $\sigma P_w$ equals  $\sigma L_w$.   
\label{prop:p247} \end{prop}

\proof (cf. \cite{Lu2}, \cite{Soibel}.) 
	Immediately follows from the discussion above.
\qed

\

\noindent More explicitly, one can follow \cite{Lu2} to identify $L_w$ with the
$\VV_w$ - orbit of $P_w$, and 
in the next section, we will define and exploit the analogues of the 
dressing action \cite{STS} to get clearer picture of the tetraplectic leaves. 
In any case,  we have the following:
 
\begin{Th} 
	The tetraplectic leaf decomposition of the quaternionic 
	flag manifold $F_n\simeq \BB\backslash\GL{\H}$ arising from the Bruhat
	quatrisson structure coincides with the Bruhat cell decomposition.  
\end{Th} 
\proof 
	The important point is that any tetraplectic leaf in $\Sp$ under 
	the quotient map $\Sp\to \Sph\backslash \Sp$ maps tetraplectomorphically 
	onto a Bruhat cell as follows from the results in this section. 
\qed

\section{Quatrisson action and intrinsic derivative}

We elaborate on some general notions related to group actions in the
quatrisson context where we recall the notation set forth
in Definition \ref{def:mult-act} and assume that we have a  multiplicative action.
Denoting   the Lie algebra of $H$ by $\hh$, we 
let:
\begin{displaymath} 
	\gamma: \hh\to \Chi(P)
\end{displaymath}
be the usual Lie algebra  anti-homomorphism, and 
recall the \emph{intrinsic derivative} of $\xi$ at $e$:
\begin{displaymath}
	d_e\xi: \hh\to \Wedge^4 \hh.
\end{displaymath}
We also define the $4$-bracket $[\cdot, \cdot, \cdot, \cdot]$ 
on $\hh^*$ to be the dual of $d\xi_e$. 
The next statement and its proof are analogous to Theorem 2.6 of \cite{LW}. 

\begin{Th} 
	In the above situation 
	for each $X\in\hh$ we have:
	\begin{displaymath}
		\LL_{\gamma(X)}\xi = \wedge^4\gamma(d_e\mu)(X).
	\end{displaymath}
	Moreover, for any $1$-forms $\omega_i$ for $1\le i \le 4$ on $P$ we have:
	\begin{displaymath}
		\LL_{\gamma(X)}\xi(\omega_1\wedge \omega_2 \wedge \omega_3 \wedge \omega_4) = 
		<[\zeta_1,\zeta_2,\zeta_3,\zeta_4],X>
	\end{displaymath}
	where $\zeta_i$ is the $\hh^\ast$-valued function on P defined by:
	\begin{displaymath}
		<\zeta_i,X> = <\omega_i,\gamma(X)> 
		\ \   \mbox{for} X \in \hh 
	\end{displaymath}
	and $[\zeta_1,\zeta_2,\zeta_3,\zeta_4]$ 
	denotes the point-wise $4$-bracket in $\hh^\ast$. 
	\label{th004}
\end{Th}

\section{Dressing action}

The Iwasawa  decomposition of $\GL{\H}=\RR\UU \Sp = \Sp \RR \UU $ 
allows us to define:
\begin{definition} 
	The \emph{dressing action} of $\RR\UU$ on $\Sp$ is the map
	$\RR\UU\times \Sp \rightarrow \Sp $ given by
	$(G,K)\mapsto K'$ where $G \cdot K = K' \cdot R \cdot  U$ for the unique 
	$R\in\RR$ and $U\in\UU$.
\end{definition}

Our goal of this section is to relate the orbits of
the dressing action with the tetraplectic leaves of the group $\Sp$.
Notice that we have restricted the usual dressing action to $\RR\UU$ since 
we will be only concerned with $\RR\UU$ orbits of the dressing
action in the remainder.
Finally, we can state the main result of this section. 
 
\begin{Th}
	The tetraplectic leaves  of  $\sham$ on $\Sp$ are the orbits of the 
	dressing action of $\RR\UU$ on $\Sp$.
\end{Th} 
\proof 
	We already know that the leaves are parametrized by $S_n$ and $\Sph$. More
	precisely, we define the center of any leaf as the element $P_w\sigma$, where $w$ as
	usual is the permutation matrix  corresponding to $w\in S_n$, and
	$\sigma\in\Sph$. The dressing action
	can be rewritten as $(G,K)\mapsto GKG'\in \Sp$, for $G,G'\in\RR\UU$, which leaves
	us in the same (open) submanifold of the Bruhat decomposition. 

	Taking $K=\sigma P_w$, we see
	that the orbit of a dressing action on a cell remains in that cell as
	there are no permutations appearing in $G$ or $G'$.  Further, the fact that the
	orbit is contained within a single leaf follows from $G$ and $G'$  being
	upper triangular with real diagonal, and thus the dressing action does not introduce
	any non-trivial elements of $\Sph$.

	For the opposite inclusion, suppose we are given two  points, $K_1,K_2$,
	of a tetraplectic leaf. As the $K_i$ are in the
	same leaf, this implies that the $K_i$'s have the same permutation 
	type, $w$,  in the Bruhat decomposition, so we write $K_i=B_i P_w V_i$ for
	some $B_i \in \BB$ and some $V_i \in\VV_w$.  Then we have
	$B_2 B_1^{-1} K_1 = K_2 V_2 V_1^{-1}$ with  $V_2 V_1^{-1} \in \VV_w$. 
	Now, as $B_i \in \BB$, we may write:
	\begin{displaymath}
		B_i = \diag{d_1^i,\dots, d_n^i}\diag{r_1^i,\dots,r_n^i}, \ \ \ r_j^i\in\R_+,
	 	\ \ d_j^i\in \Sp[1].
	\end{displaymath}
	But as the orbits are parametrized by $\Sph$, we know that $d_j^1$  
	corresponds to $d_j^2$, which
	implies that $B_2 B_1^{-1}$ must be in $\RR$ from which it follows that 
	the $K_i$'s lie in the same orbit. 
\qed

\

Another possible proof of the above result can be obtained using the
infinitesimal computations near the centers of each leaf \cite{LW},
\cite{STS}. Once we know that the tetraplectic leaves go along the orbits 
of the dressing action infinitesimally, the analyticity of the manifolds in 
question will provide a global coincidence. 

We have established that the orbits of the dressing action of $\RR\UU$ on $\Sp$ 
coincide with the tetraplectic leaves induced by the 4-vector field
$\sham$, and these are permuted by the action of $\Sph$. Therefore 
we have obtained a geometric orbit picture for any tetraplectic leaf or a 
Bruhat cell, in $F_n$.

\section{Further remarks.} 

First of all, the approach that we pursued in the present paper can be easily 
extended to all partial quaternionic flag manifolds, in particular the Grassmannians 
and projective spaces. 

It would be interesting to express the dressing action as a quatrisson action, 
with respect to a multiplicative $4$-vector field on 
$\RR\UU$. While it is clear that such a
structure exists, it is not easy to write down a local expression. It seems 
plausible that a suitable generalization of Lu-Ratiu construction \cite{LuRat}
would help. 

Evens and Lu \cite{EL} showed that the Kostant harmonic forms \cite{Kost} on complex 
flag manifolds have a Poisson harmonic nature with respect to the Bruhat Poisson
structure. It would be interesting to see how their ideas can be applied to our
situation. One can use the operator 
$\partial_\Sham= - d\circ \iota_\Sham + \iota_\Sham\circ d +\iota_\sigma$ 
to define \Sp-harmonic forms on the quaternionic flag manifolds.
Here $\sigma$ is the modular tri-vector field given by 
$d(\iota_\Sham\psi^m)=\iota_\sigma\psi^m$, and $\psi^m$ is a $\Sp$-invariant 
volume form on $F_n$. Analogously to the $T$-equivariant cohomology of complex
flag manifolds, one can consider the $\Sph$-equivariant cohomology. Another
possibility is to consider quaternionic flag manifolds as fixed
point sets of certain natural involutions on complex partial flag manifolds,
where the dimensions of the subspaces are even, and restrict certain subalgebra
of forms.  

Another possible venue to pursue is to study the hamiltonian type dynamics associated 
with the quatrisson structures. In particular, it seems that
to determine a system subject to a $\Sph$-action  which preserves a hamiltonian, 
we may need fewer integrals than in the standard Poisson case. We suspect that certain 
symmetric spaces such as quaternionic Grassmannians will have the 
property that an invariant quatrisson structure is compatible with the
Bruhat quatrisson structure, i.e. $[\chi, \Sham]=0$. This would lead to a
generalized bi-hamiltonian type systems, which are worth investigating. 

The 4-bracket on $\hh^*$ that we briefly mentioned in Section 6, gives rise 
to a certain deformed algebra of functions on $H$ (by way of the Kontsevich
formality theorem) where the deformation parameter $\hbar$ now has degree $2$.  
This implies that 
the $m_2$ term in the operadic expansion is just the standard multiplication, 
$m_3$ is trivial, and $m_4$ is determined by the bracket. This is the first 
natural occurrence of the generalized quantum group setup that we are aware of, 
and thus it seems plausible that it would lead to new interesting algebraic
structures.

\vskip 0.3in 
Department of Mathematics \\ University of Arizona\\ Tucson, AZ 85721-0089 
\\ foth@math.arizona.edu \\ litlfred@math.arizona.edu


\thebibliography{123} 

\bibitem{Perelomov}{J. A. de Asc\'arraga, A. M. Perelomov, and 
J. C. P\'erez Bueno. The Schouten-Nijenhuis bracket, cohomology
and generalized Poisson structures. {\it J. Phys. A}, {\bf 29}: 
7993-8009, 1996.}

\bibitem{Dieu}{J. Dieudonn\'e. Les d\'eterminants sur un corps non
commutatif. In French. {\it Bull. Soc. Math. France}, 
{\bf 71}: 27-45, 1943.}

\bibitem{Draxl}{P. K. Draxl. Skew Fields. {\it London Math. Soc. Lect. Not. Ser.},
{\bf 81}, Cambridge University Press, 1993.} 

\bibitem{Dr}{V. Drinfeld. Hamiltonian structures on Lie groups, Lie bialgebras
and the geometric meaning of the classical Yang-Baxter equation. 
{\it Soviet Math. Dokl.}, {\bf 27}: 68-71, 1983.}

\bibitem{EL}{S. Evens and J.-H. Lu. Poisson harmonic forms, Kostant harmonic
forms, and the $S^1$-equivariant cohomology of $K/T$. {\it Advances in Math.}, 
{\bf 142}: 171-220, 1999.} 

\bibitem{F1}{P. Foth. Tetraplectic structures, tri-moment maps, and quaternionic 
flag manifolds. {\it J. Geom. Phys.}, {\bf 41}: 330-343, 2002.}  

\bibitem{Ibanez2}{R. Ib\'a\~nez, M. de Le\'on, J. C Marrero, and E. Padr\'on. 
Nambu-Jacobi and generalized Jacobi manifolds. {\it J. Phys. A: Math. Gen.}, 
{\bf 31}: 1267-1286, 1998.}

\bibitem{KS}{L. Korogodski and I. Soibelman. Algebras of functions on quantum 
groups. {\it Math. surveys and monographs}, {\bf 56}, AMS, 1998.} 

\bibitem{Kost}{B. Kostant. Lie algebra cohomology and generalized Shubert cells.
{\it Ann. Math.}, {\bf 77}: 72-144, 1963.}

\bibitem{LuRat}{J.-H. Lu and T. Ratiu. On the non-linear convexity theorem of
Kostant. {\it Journal of AMS}, {\bf 4}: 349-363, 1991.}

\bibitem{LW}{J.-H. Lu and A. Weinstein. Poisson-Lie groups, 
dressing transformations, and Bruhat decompositions. {\it J. Diff. Geom.}, 
{\bf 31}, 501-526, 1990.}

\bibitem{Lu2}{J.-H. Lu. Coordinates on Shubert cells, Kostant's harmonic forms,and 
the Bruhat Poisson structures on $G/B$. {\it Transform. groups}, 
{\bf 4}, 355-374, 1999.}

\bibitem{Mass}{W. S. Massey. Non-existence of almost complex structures on
quaternionic projective spaces. {\it Pacific J. Math.}, {\bf 12}: 1379-1384, 1962.} 

\bibitem{STS}{M. A. Semenov-Tian-Shansky. Dressing transformations and Poisson Lie
group actions. {\it Publ. RIMS}, {\bf 21}: 1237-1260, Kyoto University, 1985.}

\bibitem{Soibel}{Y. Soibelman. 
The algebra of functions on a compact quantum group 
and its representations. {\it Leningrad J. Math.}, {\bf 2}: 161-178, 1991.}

\bibitem{Takh}{L. Takhtajan. On foundation of the generalized Nambu mechanics. 
{\it Comm. Math. Phys.}, {\bf 160}: 295-315, 1994.} 

\bibitem{Vais}{I. Vaisman. Lectures on the geometry of Poisson manifolds. 
{\it Progress in Mathematics}, {\bf 118}, Birhauser, Boston, 1994.}
\end{document}